\def\sqr#1#2{{\vcenter{\hrule height.#2pt
      \hbox{\vrule width.#2pt height#1pt \kern#1pt
         \vrule width.#2pt}
       \hrule height.#2pt}}}
\def\eop{\mathchoice\sqr34\sqr34\sqr{2.1}3\sqr{1.5}3}

\def\N{\hbox{\rm I\hskip -0.14em N}}

\def\mod{\hbox{\rm mod}\hskip 0.5em}
\def\th #1 #2. #3\par{\medbreak{\bf#1 #2.
\enspace}{\sl#3}\par\medbreak}
\documentclass{llncs1}
\usepackage{setspace}
\usepackage{makeidx}  
\usepackage{graphicx}
\begin{document}
\frontmatter          
\pagestyle{headings}  

\title{A Theorem of Congruent Primes}
\titlerunning{Of Congruent Primes}  
%
\author{Jorma Jormakka ,Sourangshu Ghosh}
\authorrunning{Jorma Jormakka ,Sourangshu Ghosh}   
%
\tocauthor{Jorma Jormakka ,Sourangshu Ghosh}

\institute{ \email{}}

\maketitle              

\begin{abstract}
 To determine whether a number is congruent or not is an old and difficult topic and progress is slow.The paper presents a new theorem when a prime number is congruent number or not. The proof is not necessarily any simpler or shorter than existing proofs, but the method may be useful in other contexts. The proof of Theorem 1 tracks the set of solutions and this set branches as a binary tree. Conditions set to the theorem restricts the branches so that only one branch is left. Following this branch gives either a solution or a contradiction. In Theorem 1 it leads to a contradiction. The interest is in the proof method, which maybe can be generalized to non-primes. 
\end{abstract}

\begin {keywordname}
Congruent numbers, elliptic curves, number theory.
\end{keywordname}

\section{Background}

A positive integer which can be written as the area of a right triangle with three rational number sides is called a Congruent Number[1]. Alternatively it can be defined as the numbers $(a,x,y,z,t)$ such that the following condition hold:
$$x^2+ay^2=z^2$$
$$x^2-ay^2=t^2$$
From this tuple $(x, y, z, t)$, we can also derive the sides of the right angle triangle ${a, b, c}$ such that
$a^2+b^2=c^2,$ and $ab/2=n$ by substituting
$$a=(y-z)/t,b=(y-z)/t,c=2x/t$$
A problem of significant interest is to determine whether a given natural number $n$ can be the area of a right-angled triangle with rational number sides.This problem
can be alternatively said of the existence of rational points on some elliptic curves that are defined over $Q$. 

Note that if we multiply each side of the triangle whose area is a congruent number $q$ by $s$, then it is evident that $s^2q$ is also a congruent number for any natural number $s$. Therefore a residue of the number $q$ in the group $Q^*/Q^{*2}$ decides whether the number $q$ will be a congruent number or not. For this reason we only consider square-free positive congruent numbers. An easy to way to determine whether a given rational number is a congruent number is the Tunnell's theorem named after number theorist Jerrold B. Tunnell who demonstrated the method in [2].

For a given square-free integer $n$, define the following numbers. 
$$A_n=((x,y,z)\in{Z^3}|n=2x^2+y^2+32z^2),$$
$$B_n=((x,y,z)\in{Z^3}|n=2x^2+y^2+8z^2),$$
$$C_n=((x,y,z)\in{Z^3}|n=8x^2+2y^2+64z^2),$$
$$D_n=((x,y,z)\in{Z^3}|n=8x^2+2y^2+16z^2),$$

Now if $n$ is actually a congruent number then by Tunnell's theorem we have if $n$ is odd then $2A_n = B_n$ and if $n$ is even then $2C_n = D_n$. There have many works to classify congruent numbers especially for primes. Gross[3] proved that if $n$ is square free integer and has at most two
prime factors of the form 5,6 or 7 (mod 8), then $n$ must be a congruent number. Monsky[4] proved the following important properties to determine whether a number is congruent.The following are all congruent numbers:
 
\begin{enumerate}
    \item $p_5, p_7, 2p_7,$ also proved by Stephens[5]
    \item $2p_3$ also proved by Heegner[6], and Birch (1968)[7,8]
    \item $p_3p_5, p_3p_7, 2p_3p_5, 2p_5p_7$
    \item $p_1p_5$ provided $(p_1/p_5)=-1$ holds
    \item $p_1p_3$ provided $(p_1/p_3)=-1$ holds
    \item $p_1p_7,2p_1p_7$ provided $(p_1/p_7)=-1$ holds
\end{enumerate}
Here $p_k$ refers to an arbitrary prime congruent to $k$ mod 8 and $(a/p)$ (where gcd($a$,$p$)=1) is the Legendre symbol which is 1 and -1 if $a$ is a quadratic residue of $p$ or not respectively. Iskra[9] proved the following important properties to determine whether a number is not congruent. The following are not congruent numbers:
\begin{enumerate}
    \item $p_3, 2p_5, p_3q_3, 2p_5q_5$ also proved by Genocchi[10].
    \item $p_3^1p_3^2..p_3^t$ provided $(p_3^m/p_3^n)=-1$ holds for $m<n$.
\end{enumerate}
Lagrange[11] similarly proved the following important properties for a number to be a non-congruent number
\begin{enumerate}
\item $p_1p_3$ provided $(p_1/p_3)= -1$ holds 
\item $2p_1p_5$ provided $(p_1/p_5)= -1$ holds
\item $n = p_1p_3q_1$ with the condition that $n$ can be written as $n = pqr$ or $2pqr$ such that $(p/q)=(p/r)=-1$
\item $n = 2p_1p_5q_1$ with the condition that $n$ can be written as $n = pqr$ or $2pqr$ such that $(p/q)=(p/r)=-1$
\end{enumerate}
Bastien[12] similarly proved the following important properties for a number to be a non-congruent number 
\begin{enumerate}
     \item if $p\equiv 9 \ (\mod 16)\$$ then $2p$ is a congruent number.
     \item $n=p_1,p_1=a^2+4b^2, ((a+2b)/p_1)=-1$
\end{enumerate}
Tian Ye[13] proved that for a given number $k$ in the congruence classes 5, 6, 7 (mod 8),there are infinitely many square-free congruent numbers with $k$ prime factors.

\section{The New Method}

Consider an elliptic curve of the form:
$$y^2=x^3-d^2x\eqno(1)$$
where $d$ is an integer. A rational solution $(x,y)$ to the elliptic
curve (1) is a solution where $x$ and $y$ are rational numbers.   

The substitution 
$x=d(a+b)/b$, $y=2d^2(a+c)/b^2$ 
changes $y^2=x^3-d^2x$ to 
$a^2+b^2=c^2$ with $ab=2d$. Then $4d^2=a^2(c^2-a^2)$. 
Integers $d$ that give rational number solutions 
to $a^2+b^2=c^2$, $ab=2d$ are called congruent numbers. If $d$ 
is a congruent number the elliptic curve (1) has a rational solution
where $y$ is not zero. In that case it has infinitely many rational 
solutions.

If there is a solution for $d=s^2$, then there is a solution for
$d=1$ because the substitution $y=s^3y'$, $x=s^2x'$ changes $y^2=x^3-d^2x$ 
to $y'^2=x'^3-x'$. It is known that no $d=s^2$ is a congruent number as proved by Fermat using his method of infinte descent. 
The case where $d$ is a prime number is almost solved.

For notations the following concepts suffice:
The condition that the integer $a$ divides integer $b$
is written as $a|b$. If $p>2$ is a prime, 
the cyclic group of integers modulo $p$ is denoted by $Z_p$ and 
$Z_p^*=\{1,\dots ,p-1\}$.
The set of quadratic residues modulo $p$ is the set
$$QR_p=\{ x\in Z_p^* | \exists y\in Z_p^* \ {\rm such} \ {\rm that}\\
 \ y^2\equiv x \ (\mod p)\}.$$
The set of quadratic nonresidues modulo $p$ is the set
$$QNR_p=\{ x\in Z_p^* | x\not\in QR_p\}.$$

Let us start by two very simple lemmas.

\th Lemma 1. Let $c^2=a^2+b^2$, $a,b,c\in Z$, then $\exists h,m,e\in \N$
such that
$$a=\pm hem \ , \ b= \pm {1\over 2}h(m^2-e^2) \ ,\\
\ c=\pm {1\over 2}h(m^2+e^2).$$

\proof 
Without loss of generality we can assume that $a,b,c\in N$. 
We can write $c^2-b^2=(c-b)(c+b)=a^2$. Let $h=gcd(c+b,c-b)$. Then there exists
$m$ and $e$, $m>e$, $gcd(m,e)=1$, such that
$c+b=hm^2$, $c-b=he^2$. The claim follows. $\eop$

With Lemma 1 we can characterize congruent numbers.

\th Lemma 2. Let $d\in Z$, $d>0$. 
Rational solutions $(x,y)$ with $x\not=0, y\not=0$ 
to
$$y^2=x^3-d^2x$$
are of the form
$$(x_1,y_1)=\left(\\
d{m+e\over m-e},\pm {k\over j}d{m+e\over m-e}\right),$$
$$(x_2,y_2)=\left(\\
d{m-e\over m+e},\pm {k\over j}d{m-e\over m+e}\right),$$
where $k,j,e,m\in N$, $m>e$, $gcd(m,e)=1$, $gcd(k,j)=1$, satisfy
$$d=\left({k\over 2j}\right)^2{m^2-e^2\over em}.\eqno(2)$$

\proof
Let $x,y\in Q$, $x\not=0, y\not=0$. 
Let us write $\alpha={d\over x}+1\in Q$, $\beta={y\over x}\in Q$.
Solving (10) for $x$ and solving $x$ from the definition of $\alpha$ yields
$$x={\beta^2\over 2\alpha-\alpha^2}={d \over \alpha-1}.$$
Writing $\beta={k\over j}$ for some $k,j\in N$ gives
$$\alpha_{1,2}=1-{k^2\over j^22d}\pm {\sqrt{(2dj^2)^2+(k^2)^2}\\
\over j^2 2d}.$$
As $y\not=0$, $k\not=0$.
By Lemma 1, $\alpha_{1,2}\in Q$ if and only if there exist
$h,e,m\in N$, $gcd(e,m)=1$, $m>e$, such that
$$k^2=hem \ , \ 2dj^2={1\over 2}h(m^2-e^2) \ , \ c={1\over 2}h(m^2+e^2).$$
If $em=0$, then $k=0$ and $y=0$. This solution gives $j=2dj^2$
$$\alpha_{1,2}=1\pm {2dj^2\over 2dj^2}=1\pm 1 \ , \ \alpha_1=2, \alpha_2=0,$$
$$x_1={d\over \alpha-1}=d \ , \ x_2=-d \ , y=0$$
but we have excluded this case in the assumptions.
Since $em\not=0$, let us write $h={k^2\over em}$. Eliminating $h$ yields
$$d=\left({k\over 2j}\right)^2{m^2-e^2\over em},$$
$$c={k^2\over 2}({m^2+e^2}).$$
Simplifying $\alpha_{1,2}$ yields
$$\alpha_{1,2}={1\over m^2-e^2}\left(m^2-e^2-2em\pm (m^2+e^2)\right),$$
i.e.,
$$\alpha_1={2m\over m+e} \ , \ \alpha_2=-{2e\over m-e}$$
$$x_1={d\over \alpha_1-1}=d{m+e\over m-e} \ , \ x_2=-d{m-e\over m+e},$$
$$y=\beta x \ , \beta^2=\left({k\over j}\right)^2=4d{em\over m^2-e^2}.$$
This gives the claim. $\eop$

As two examples of Lemma 2 
$$d=5=\left({3\over 2\cdot 2}\right)^2{9^2-1^2\over 9\cdot 1}$$
$$d=7=\left({24\over 2\cdot 5}\right)^2{16^2-9^2\over 16\cdot 9}$$
are both congruent numbers. Notice that $gcd(k,j)=1$ but it is allowed
that $2|k$. 

If $d$ is a square, there are no rational
solutions to (1) with $y\not=0$.
There are the three solutions $(0,0), (d,0), (-d,0)$ to (1), so the number of
rational solutions of (1) is finite, the rank of the elliptic curve is zero. 

In the next theorem gives a set of values where $d$ is a prime number 
and (1) has no rational solutions, i.e., the elliptic curve 
has rank zero. The case of prime numbers $d$ is rather well known:
if $p\equiv 5 \ (\mod d)$ or $p\equiv 7 \ (\mod d)$ 
the number $d$ is
a congruent number and there are solutions to (1). If
$p\equiv 3 \ (\mod d)$ there are no solutions and $d$ is not a congruent
number. The only case remaining is $p\equiv 1 \ (\mod d)$. 
For that case
it is known that e.g. $p=41$ is a congruent number, while e.g. $p=17$ is not.

The next theorem does not solve the problem for any prime $p$ that
is one modulo eight because if $p\equiv 1 \ (\mod 8)$ it is 
necessarily true that $-1\in QR_p$, i.e., $-1\in QR_p$ is equivalent
with the condition that $4|(p-1)$ and if $p\equiv 1 \ (\mod d)$, then
$8|(p-1)$. The theorem does prove e.g. that $p=19$ is not a congruent
number, but as $19\equiv 3 \ (\mod d)$ this is known. Yet, the method of this
proof seemed interesting enough to me in order to be written down. 
The method may generalize to other numbers than primes. The primality
condition is used only in a few places. The main idea is to exclude
branches from a recursion.

\th Theorem 1. Let $d>3$ be a prime such that $-1\in QNR_d$ and $2\in QRN_d$.
The equation (2) in Lemma 2
does not have solutions $k,j,m,e\in \N$ where $gcd(m,e)=1$, 
$gcd(k,j)=1$, $m>e>0$.

\proof
We write (2) with $m_1$, $e_1$ 
$$d=\left({k\over 2j}\right)^2{m_1^2-e_1^2\over e_1m_1}\eqno(3)$$
If $d|m_1$ then $d|e_1$ and $gcd(m_1,e_1)\not=1$, thus $d\not| m_1$ 
and $d\not| e_1$.
If $d|k^2$ then since $d$ is a prime $d|k$. It follows that $k=dk_1$
and as $gcd(k,2j)=1$ holds $d\not| 2j$. Thus
$$(2j)^2m_1e_1=dk_1^2(m_1^2-e_1^2)$$
which is not possible as the left side is not divisible by $d$. Thus
$d\not| k^2$. Therefore $d|m_1^2-e_1^2$. 

If $2\not|k$ we convert (3) into the form 
$$d=\left({k\over j}\right)^2{st\over m^2-e^2}\eqno(4)$$
by the substitution $m_1=m+e$, $e_1=m-e$, i.e., 
$2m=m_1+e_1$, $2e=m_1-e_1$. 
As $m_1e_1=(m+e)(m-e)=m^2-e^2$ holds $em={1\over 4}(m_1+e_1)(m_1-e_1)$.
As $4|(m_1^2-e_1^2)$ in (3) if $2\not|k$
it follows that one of $m_1+e_1$ or $m_1-e_1$
is even. If so, they are both even and $2|m_1+e_1$, $2|m_1-e_1$ and
$m,e$ are integers. As $gcd(m_1,e_1)=1$, $gcd(m_1+e_1,m_1-e_1)=2$.
Then $gcd(m,e)=gcd(((m_1+e_1)/2)((m_1-e_1)/2))=1$. 
Since $m_1>e_1>0$ holds $m>e>0$.

If $2|k$ then the substitution is 
$m=m_1+e_1$, $e=m_1-e_1$. Then $m,e$ are integers and $m>e>0$.
In this case $2\not|j$ gecause $gcd(k,j)=1$. Therefore $2\not| (m_1^2-e_1^2)$.
It follows that  $gdc(m,e)=gcd(m_1+e_1,m_1-e_1)=1$. We get the same form
(4) since $me=m_1^2-e_1^2$ and $m^2-e^2=4m_1e_1$.  

Then $d|em$ and $j^2|em$. Let us write (4) as
$$j^2(m+e)(m-e)d=k^2me.\eqno(5)$$
Since $gcd(m,e)=1$ it follows that $gcd(m\pm e,m)=1$. Indeed, if
$m\pm e=c_1r$, $m=c_2r$ for some $r, c_1,c_2\in \N$, then 
$$c_1c_2r=c_2m\pm c_2e=c_1m \Rightarrow (c_1-c_2)m=\pm c_2e$$
$$\Rightarrow m|c_2 \Rightarrow \exists \alpha \in \N \\
\ {\rm such \ that \ } c_2=\alpha m$$
$$\Rightarrow m=\alpha mr \ \Rightarrow \alpha r=1 \Rightarrow r=1.$$
Similarly, $gcd(m\pm e, e)=1$. 

Since $gcd(k,j)=1$ it follows from (4) that $k^2=m^2-e^2$. 
Therefore (4) implies that $dj^2=em$.
As $dj^2=em$ and $gcd(e,m)=1$ there is one of the cases: either
$m=ds^2$, $e=t^2$ for some $s,t>0$ or $m=s^2$, $e=dt^2$.

As $k^2=(m+e)(m-e)$ and $gcd((m+e)(m-e))\le 2$
we have two cases cases: either 
$m+e=c_1^2$ and $m-e=c_2^2$ 
for some $c_1,c_2>0$ or $m+e=2c_1^2$ and $m-e=2c_2^2$. 
  
We have four cases in total.

Case 1. $m=ds^2$, $e=t^2$, $m+e=c_1^2$, $m-e=c_2^2$. Then
$$m-e=s^2d-t^2=c_2^2.$$
The equation yields $-1\equiv (c_2t^{-1})^2 \ (\mod d)$
which is impossible since $-1\in QNR_d$.

Case 2. $m=ds^2$, $e=t^2$, $m+e=2c_1^2$, $m-e=2c_2^2$. Then
$$s^2d+t^2=2c_1^2 \ , \ s^2d-t^2=2c_2^2.$$
Multiplying the modular equations
$$t^2\equiv 2c_1^2 \ (\mod d) \ , -t^2\equiv 2c_2^2 \ (\mod d)$$
yields $-1\equiv (2c_1c_2t^{-2})^2 \ (\mod d)$
which is impossible since $-1\in QNR_d$.

Case 3. $m=s^2$, $e=dt^2$, $m+e=c_1^2$, $m-e=c_2^2$. Then
$$s^2+t^2d=c_1^2 \ , \ s^2-t^2d=c_2^2.$$
Thus
$$2s^2=c_1^2+c_2^2\eqno(6)$$
so
$$4s^2=c_1^2+2c_1c_2+c_2^2+c_1^2-2c_1c_2+c_2^2$$
$$(2s)^2=(c_1+c_2)^2+(c_1-c_2)^2.\eqno(7)$$
It follows from Lemma 1 that 
$\exists h', e', m' \in \N$, $gcd(m',e')=1$ such that
$$c_1+c_2=h'e'm' \ , c_1-c_2={1\over 2}h' (m'^2-e'^2),$$
$$2s={1\over 2}h' (m'^2+e'^2).$$
Solving $c_1,c_2,s$ yields
$$c_1={1\over 4}h' (2e'm'+m'^2-e'^2),$$
$$c_2={1\over 4}h' (2e'm'+e'^2-m'^2),$$
$$s={1\over 4}h' (m'^2+e'^2).$$

Since
$$2t^2d=c_1^2-c_2^2=(c_1-c_2)(c_1+c_2)$$
we get
$$d={1\over 4t^2}h'^2 e'm'(m'^2-e'^2)$$
i.e.
$$d=\left({h'e'm'\over 2t}\right)^2 {(m'^2-e'^2)\over e'm'}.$$
Removing the greatest common divisor of $h'e'm'$ and $t$ 
this equation can be written as 
$$d=\left({k_{i+1}\over 2j_{i+1}}\right)^2 {(m_{i+1}^2-e_{i+1}^2)\over e_{i+1}m_{i+1}}.\eqno(8)$$
As $gcd(m',e')=1$ and we made $gcd(k,j)=1$, equation (8) is
is of the same form as (3)
$$d=\left({k_i\over 2j_i}\right)^2 {(m_i^2-e_i^2)\over e_im_i}=\left({k\over 2j}\right)^2 {(m_1^2-e_1^2)\over e_1m_1}.$$
We have a recursion that in each step reduces the numbers $m_i,e_i$ 
to numbers $m_{i+1},e_{i+1}$ that are of the order of square root
of $m_i,e_i$.

Case 4. $m=s^2$, $e=dt^2$, $m+e=2c_1^2$, $m-e=2c_2^2$. 
We can select $c_1>c_2\ge 0$.
Then
$$s^2+t^2d=2c_1^2 \ , \ s^2-t^2d=2c_2^2.$$
Thus
$$s^2=c_1^2+c_2^2 \,\ dt^2=c_1^2-c_2^2=(c_1-c_2)(c_1+c_2).\eqno(9)$$
Let us notice that $m+e=2c_1^2$ and 
$$1=gcd(m+e,e)=gcd(2c_1^2,dt^2) \Rightarrow gcd(c_1,t)=1, gcd(2,t)=1$$ 
$$1=gcd(m-e,e)=gcd(2c_2^2,dt^2) \Rightarrow gcd(c_2,t)=1.$$

First we exclude one case in the second equation of (9).
If $t>1$ and $c_1+c_2=\alpha_1t$ and $c_1-c_2=\alpha_2 t$ for some 
$\alpha_1, \alpha_2 \in \N$, then
$$2c_1=(\alpha_1+\alpha_2)t \Rightarrow t=1, 2c_1=\alpha_1+\alpha_2,$$ 
$$2c_2=(\alpha_1-\alpha_2)t \Rightarrow t=1, 2c_1=\alpha_1-\alpha_2.$$
Thus, $dt^2=c_1^2-c_2^2=\alpha_1\alpha_2t^2$. It follows that
$d=\alpha_1\alpha_2$ and as $d$ is prime and necessarily 
$\alpha_1>\alpha_2$ it follows that $\alpha_1=d$, $\alpha_2=1$.
Then $c_1=d+1$ and $c_2=d-1$. Consequently $s^2=c_1^2+c_2^2=2(d^2-1)$
is even, so $m$ is even. Since $s^2+dt^2=2c_1^2$ it would follow that
$t$ is also even as $d$ is odd, but $t=1$ in this case. We have a 
contradiction. 

Thus, in (9) must be one of the three cases
$$t^2 | (c_1+c_2) \Rightarrow (c_1-c_2)|d \\
\Rightarrow c_1-c_2=d\Rightarrow t^2=c_1+c_2,$$
or 
$$t^2 | (c_1-c_2) \Rightarrow (c_1+c_2)|d \\
\Rightarrow c_1+c_2=d\Rightarrow t^2=c_1-c_2,$$
or
$$t=1.$$
In the first case 
$$2c_1=t^2+d\ge 0 \ , \ 2c_2=t^2-d\ge 0.$$
In the second case
$$2c_1=d+t^2\ge 0 \ , \ 2c_2=d-t^2\ge 0.$$
In both of these two cases we can derive in a similar way:
$$s^2=c_1^2+c_2^2 \ \Rightarrow (2s)^2=(2c_1)^2+(2c_2)^2$$
yields
$$(2s)^2=(d+t^2)^2+(d-t^2)^2.\eqno(10)$$
By Lemma 2 there exist 
$h',e',m'\in \N$ such that
$$d+t^2=h'e'm' \ , \ d-t^2={1\over 2}h'(m'^2-e'^2).$$
The first equation implies that $d\not| h'$. 
Thus
$$4d=h'((m'+e')^2-2e'^2)$$
i.e., as $h'\not\equiv 0 \ (\mod d)$
$$2\equiv (m'^2+e'^2)^2e'^{-2} \ (\mod d)\eqno(11)$$
which is a contradiction since $2\in QNR_d$.
There remains the case $t=1$. Then $2c_1^2=s^2+d$, $2c_2^2=s^2-d$.
Instead of (10) we get
$$(2s)^2=(d+s^2)^2+(d-s^2)^2.$$
The contradiction (11) comes in the same way with $t$ replaced by $s$.
This means that Case 4 is not possible.    

Because Cases 1, 2 and 4 are not possible, only Case 3 is left. Case 3
gives a recursion formula.
The values $h',m',e'$ in Lemma 1 satisfy
$${e'\over m'}={a\over b+c}={c-b\over a}$$
$$h'=gcd(b+c,b-c)$$
giving $a^2=c^2-b^2$. The numbers $h',m',e'$ can be chosen to be 
positive and on the order of $a,b,c$. Thus, $h',m',e'$ in (8) are of the order
$c_1,c_2$. The numbers $c_1,c_2$ are of the order $\sqrt{m}$, $\sqrt{e}$.
Therefore in each step the numbers $m_i,e_i$ get smaller, they are reduced 
to the order of their square roots. Consider the problem when the recursion
stops.

Let us look at an example of $d=5$. Then 
$$d=5=\left({3\over 2\cdot 2}\right)^2{9^2-1^2\over 9\cdot 1}.$$
We have $m_1=9,e_1=1,k=3,j=2$. 
We can do the first step and find $m=5,e=4$ and 
$$d=5=\left({3\over 2}\right)^2{5\cdot 4\over 5^2-4^2}.$$
Identifying $k^2=3^2=5^2-4^2=9$, $j^2d=4\cdot 5=20=5\cdot 4=me$, 
$m=ds^2=5\cdot 1^2$,
$e=t^2=2^2$, $m+e=5+4=3^2=c_1^2$ and $m-e=5-4=1^2=c_2^2$ shows that
the logic in the lemma is correct. We have Case 1, but for $d=5$ the
conditions of the lemma are not fulfilled:
$-1\in QNR_5$. This is why Case 1 does not give a contradiction. What happens
in Case 1 is that when we remove the term $dt^2$ in a case resembling (6)
we do not get (6) but 
$$2t^2=c_1^2-c_2^2$$
Therefore we do not get (7) which can be inserted to 
the equation to Lemma 1 for calculation of the numbers $h',m',e'$. 

Let us look at another example, that of $d=7$. Here $-1\in QNR_7$ and
the Case is not 1. 
$$d=7=\left({24\over 2\cdot 5}\right)^2{16^2-9^2\over 16\cdot 9}.$$
We have $m_1=16,e_1=9,k=24,j=5$. We find $m=16+9=25,e=16-9=7$. Thus 
$$d=7=\left({24\over 5}\right)^2{25\cdot 7\over 25^2-7^2}.$$

Here $k^2=24^2=576=25^2-7^2=m^2-e^2$, $j^2d=25\cdot 7=175=25\cdot 7=me$, 
$m=s^2=5^2$, $e=dt^2=7\cdot 1^2$, 
$m+e=25+7=32=2\cdot 4^2=2c_1^2$ and $m-e=25-7=18=2\cdot 3^2=2c_2^2$. The
Case is 4. We notice that $t^2=1$ and $c_1=4$, $c_2=3$, thus we have
the case $t=1$. Then $s^2+d=5^2+7=32=2\cdot 4^2=2c_1^2$ and
$s^2-d=5^2-7=18=2\cdot 3^2=2c_2^2$. We get
$$(2s)^2=100=64+36=(2c_1)^2+(2c_2)^2=(5^2+7)^2+(5^2-7)^2$$
and therefore find the numbers $h',m',e'$ for $10^2=8^2+6^2$.
The numbers are $h'=gcd(10+6,10-6)=4$, $e'=1$, $m'=2$. Thus
$$d+t^2=h'e'm'=7+1=8 \ , \ d-t^2={1\over 2}h'(m'^2-e'^2)=6$$
are true and
$$4d=h'((m'+e')^2-2e'^2)=28=4\cdot (3^2-2).$$ 
We get the modular equation $3^2\equiv 2 \ \mod(7)$, which
violates the assumption $2\in QNR_d$, but indeed $2\in QR_7$. 
Therefore for $d=7$ we do not get a contradiction. 

The way the lemma works is that in (2)
the numbers $m_1$ and $e_1$ must be squares $m_1=s_1^2$, $e_1=t_1^2$ 
so that $k^2$ can cancel them. The condition $-1\in QNR_d$ excludes
the larger branch $(s_1^2+t_1^2)$ of 
$$m_1^2-e_1^2=(s_1^2+t_1^2)(s_1^2-t_1^2)$$
by $(s_1^2+t_1^2)\equiv 0\hskip 1em (\mod d)$ being impossible. 

Therefore
$4d|(m_1^2-e_1^2)$ leads to $4d|(s_1^2-t_1^2)$.
The condition $2\in QNR_d$ excludes Case 4 and leaves only Case 3 which
gives a recursion. Thus, the numbers $m_i,e_i$ get smaller.

If there is a congruent number $d$ with $-1\in QNR_d$, the recursion must
continue until it stops in some way and not to a contradiction, 
but the recursion does not stop and continues to a contradiction. 
At each stage
$4d|(m_i^2-e_i^2)$ or $d|(m_i^2-e_i^2)$ depending on if $k_i$ is odd or even.
The numbers $m_i$ and $e_i$ become smaller on each 
step. Finally we must have $4d=m_i^2-e_i^2$ or $d=m_i^2-e_i^2$.

Changing variables in (2) to
$m=(m_i+e_i)/2$, $e=(m_i-e_i)/2$ if $k$ is odd and
$m=m_i+e_i$, $e=m_i-e_i$ if $k$ is even
we get
$$d={k^2\over j^2}{me\over m^2-e^2}.\eqno(12)$$
When the recursion has reached $4d=m_i^2-e_i^2$ or $d=m_i^2-e_i^2$
the number $j=1$. 
In (12) necessarily $k^2=m_i^2e_i^2$ and consequently $d=me$.
As $d$ is prime either $m=d$, $e=1$ or $m=1$, $e=d$. As in Cases 1 and 2
the choice $m=d$ leads to $-1\in QR_d$ and is impossible. 
Thus $m=1$ and $t=d$, but then $m^2-e^2<0$ and $d>0$ is negative. This
is a contradiction. The recursion leads to a contradiction and
the claim of the lemma follows. $\eop$

There are primes $d$ filling the conditions of the lemma: for $d=19$
holds $-1\in QNR_{19}$ and $2\in QNR_{19}$. We also get a small result:

\th Corollary 1.
If $p$ is a prime and $p\equiv 7 \ (\mod 8)$, then 
$2\in QR_p$.

\proof
If $p$ is a prime and $p\equiv 7 \ (\mod 8)$, then $p$ is a congruent
number. Therefore the conditions of Theorem 1 cannot be fulfilled.
The condition $-1\in QR_p$ is equivalent with $4 | (p-1)$. 
As $p-1=6+8k$ for some $k$, it follows that $4 \not| (p-1)$. Thus
$-1\in QNR_p$. The only other condition in Theorem 1 is that $2\in QNR_p$. 
$\eop$

This can be otherwise be proved easily without using the Theorem 1 by Gauss' lemma which states that $(a/p)=(-1)^n$. Here $(a/p)$ (where gcd($a$,$p$)=1) is the Legendre symbol which is 1 and -1 if $a$ is a quadratic residue of $p$ or not respectively. Here n is the number of integers in the set
$$S=(a,2a,3a,....,((p-1)/2)a)$$
whose remainder will be greater than $p/2$ when divided by $p$.
Putting the value of $a$ as 2 we get $(2/p)=(-1)^n$, where n is the number of integers in the set
$$S=(1,2.1,3.1,....,((p-1)/2).2)$$

Note that all of the elements present in $S$ is smaller than $p$. Therefore the problem reduces to only count the number of elements that exceed $p/2$.The number of such integers will be $n=(p-1)/2-[p/4]$. If $p$ is a prime and $p\equiv 7 \ (\mod 8)$, we have $p$ of the form $8k+7$. Therefore
$$n=(8k+7-1)/2-[8k+7/4]=4k+3-(2k+1)=2k+2$$
As $n$ is even, we have $2\in QR_p$.

\section{Conclusion}

Whether primes are congruent numbers of not is an old and difficult topic and progress is slow. The paper presents a new proof to a known theorem. The proof is not necessarily any simpler or shorter than existing proofs, but the method may be useful in other contexts. The proof of Theorem 1 tracks the set of solutions and this set branches as a binary tree. Conditions set to the theorem restricts the branches so that only one branch is left. Following this branch gives either a solution or a contradiction. In Theorem 1 it leads to a contradiction. Using different conditions in this method may give new results.


\begin{thebibliography}{1}
%


\bibitem 
{Koblitz}
Koblitz, N. Introduction to Elliptic Curves and Modular Forms. New York: Springer-Verlag, 1993.
\bibitem 
{Tunell}
Tunnell, Jerrold B. (1983), "A classical Diophantine problem and modular forms of weight 3/2", Inventiones Mathematicae, 72 (2): 323–334, doi:10.1007/BF01389327, hdl:10338.dmlcz/137483
 \bibitem
 {Gross}
 Gross, B. H., Zagier, D. B. (1986) Heegner points and derivatives of L-series, Invent.
Math, 84(2), 225–320
\bibitem
{Monksy}
 Paul Monsky (1990), "Mock Heegner Points and Congruent Numbers", Mathematische Zeitschrift, 204 (1): 45–67, doi:10.1007/BF02570859
  \bibitem
 {Stephens}
 Stephens, N. M. (1975) Congruence properties of congruent numbers, Bull. London Math. Soc., 7, 182–184.
  \bibitem
 {Heegner}
Heegner, K. (1952) Diophantine analysis und Modulfunctionen, Math. Z., 56, 227–253.
  \bibitem
 {Birch1}
 Birch, B. J. (1968) Diophantine analysis and modular functions, Proc, Bombay Colloq. Alg. Geom., 35–42.
  \bibitem
 {Birch2}
 Birch, B. J. (1970) Elliptic curves and modular functions, Symp. Math. 1st. Alta MAt., 4,
27–32.

\bibitem
 {Iskra}
 Iskra, Boris. Non-congruent numbers with arbitrarily many prime factors congruent to $3$ modulo $8$. Proc. Japan Acad. Ser. A Math. Sci. 72 (1996), no. 7, 168--169. doi:10.3792/pjaa.72.168. https://projecteuclid.org/euclid.pja/1195510284
 \bibitem
 {Genocchi}
  A. Genocchi. Note analitiche sopra tre scritti. Annali di Scienze Matematiche e Fisiche, 6:273–317, 1855. → page 5.
 \bibitem
 {Lagrange} 
  J. Lagrange, Nombres congruents et courbes elliptiques, S´em. Delange–Pisot–
Poitou, 16e ann´ee, 1974/75, no. 16.
 \bibitem
 {Bastien} 
  L.Bastien,L'intermediaire des math,21,1914,20-21,231-2
 \bibitem
{Tian ye} 
 Tian, Ye (2014), "Congruent numbers and Heegner points", Cambridge Journal of Mathematics, 2 (1): 117–161, arXiv:1210.8231, doi:10.4310/CJM.2014.v2.n1.a4, MR 3272014.




\end{thebibliography}
\end{document}